\magnification=1200
\baselineskip=14pt

\def\CMP {1}
\def \BVa {2}
\def \BVb {3}
\def\JMPa {4}
\def\MPL{5}
\def\JMPd {6}
\def\SchCMPa{7}
\def \Gayduk {8}
\def\But {9}
\def\Poin {10}
\def\BarScha{11}
\def\BarSchb{12}
\def\leites{13}
\def\SchCMPb {14}
\def\Berezin {15}
\def\leitb {16}
\def\Shander{17}
\def\buttin{18}
\def\Gui  {19}

\def\p {\partial}
\def \D {\Delta_{d{\bf v}}}
\def \Ds  {\Delta^{\#}}
\def\t{\tilde}
\def\s {\sigma}
\def\L {\Lambda} 
\def\Darboux {$z^A=$  $x^1,\dots,x^n$, $\theta_1,\dots,\theta_n$}   
\def\a{\alpha}
\def\O{\Omega}
\def\d{\delta}
\def\dv  {{d{\bf{v}}}}
\def\A {{\cal A}}

$\qquad$ $\qquad$ $\qquad$ $\qquad$
$\qquad$ $\qquad$ $\qquad$ $\qquad$
                    math/9909117

     \centerline {\bf $\Delta$-Operator on Semidensities} 
            \centerline     {\bf and} 
	\centerline {\bf Integral Invariants in 
	 the Batalin-Vilkovisky Geometry}

                      \medskip
                   \centerline {{\bf O.M.Khudaverdian}
		   \footnote
		   {$^*$}
	{e-mails:"khudian@thsun1.jinr.ru", 
"khudian@vxjinr.jinr.ru,"khudian@sc2a.unige.ch"}}
		\medskip   
\centerline {\it Laboratory of Computing Technique and Automation}
\centerline   {\it Joint Institute for Nuclear Research}
    \centerline {\it Dubna, Moscow Region  141980 Russia}        
 \centerline {\it on leave of absence of the Department of Theoretical
       Physics of Yerevan State University}
  \centerline{A. Manoukian St., 375049  Yerevan, Armenia}
        \medskip
       \centerline {September, 1999}

                \medskip

  {\bf Abstract} {\it The action of Batalin-Vilkovisky
  $\Delta$-operator on semidensities in an
  odd symplectic superspace is defined.
   This is used  for the construction
   of integral invariants on surfaces embedded
in an odd symplectic superspace 
and for more clear interpretation of the
 Batalin-Vilkovisky formalism
 geometry.} 

\bigskip
                 \centerline{\bf 1.Introduction}
		 
                  \medskip
		  
   Density of the weight $\s$ is the function on a space (superspace)
   subject to the condition that under change of coordinates it is
   multiplied on the $\s$-th power of the determinant (Berezinian)
    of the transformation. The density of the weight $\s=1$ is
    a volume form on the space.

 In this paper 
   we consider  semidensities (densities of the weight $\s=1/2$)in
  a  superspace provided with an odd symplectic structure and
  define the action of $\Delta$-operator on them.  
  Using these
 constructons we come to the new outlook on
 the invariant semidensity 
  in an odd symplectic superspace [\CMP] 
   and construct integral invariants on surfaces
   of codimension $(1.1)$ embedded in an odd sympelctic
    superspace.
    We also analyze from this point of view the
     geometrical formulation of the Batalin-Vilkovisky (BV) formalism.

  The concept of an odd symplectic superspace
  and $\Delta$-operator on it appeared in mathematical physics 
 in the pioneer works of Batalin-Vilkovisky [\BVa,\BVb], where these objects
were used for constructing covariant Lagrangian version of the
  BRST quantization. The geometrical meaning of these objects
    and interpretation of BV master-equation in its terms
    was studied in [\JMPa,\MPL,\JMPd]. 
  The complete
 analysis of these constructions and their meaning in BV formalism
  was performed by A.S.Schwarz  
in  [\SchCMPa]. In particular in this paper the
  essential role of the semidensity 
  in an odd symplectic superspace was studied. It 
   turned out to be a volume form on Lagrangian surfaces
   in this superspace.
  
      Let us shortly sketch these results.
    
    We call a superspace provided 
    with an odd symplectic structure
     {\it odd symplectic superspace}. 
      An odd symplectic superspace
     which is provided in addition with a volume form
    will be called {\it special odd symplectic superspace}.
     
     If $d{\bf v}$ is the volume form of the special odd symplectic 
     superspace, then one can consider operator 
     $\D$ whose action on a
     function on this superspace is equal  to the divergence
     w.r.t. the volume form $d{\bf v}$ of the Hamiltonian vector field
     corresponding to this  function.
    This second order differential operator is not trivial because
      transformations
      preserving odd symplectic structure do not preserve
      any volume form (Liouville theorem fails to be fulfilled in
      a case of odd symplectic structure). 
      
       We call a superspace 
       {\it normal special odd symplectic superspace}
       if in the special odd symplectic superspace 
       there exist such Darboux coordinates, 
	  that in these coordinates
        the density of the volume form is equal to one:
                  	$$
	\dv=dx^1\dots dx^n d\theta_1\dots d\theta_n\,.
	   \eqno (1.1)
	        $$
        (Darboux coordinates in an odd symplectic superspace
	are  coordinates \Darboux $\,$ in which the 
	 Poisson bracket corresponding to the symplectic structure has
	 the canonical form:
	 $\{x^i,\theta_j\}=\delta^i_j$,
	  $\{x^i,x^j\}=0$.)

	The concept of normal special odd symplectic
	superspace,
	which is called also  SP superspace ([\SchCMPa]),
	 is crucial in the geometrical interpretation
	  of BV formalism.
	
	If $f$ is an even function on 
	a normal special odd symplectic superspace
	and $d{\bf v}^\prime=f d{\bf v}$ is
	a new volume form on it then
        the main essence of geometrical
	 formulation of BV formalism
	 can be shortly expressed in the following two statements:
	
         {\bf Statement 1.}  
	 
The following conditions:	
                      $$
		     \matrix
		        { 
  a)\quad	\hbox{  the volume form $d{\bf v}^\prime$ provides}\cr 
   \hbox{the odd symplectic superspace by the normal special 
   structure}\cr
   \hbox{ in the same way as the initial volume form $\dv$}, \cr
   \hbox{(i.e. there exist Darboux coordinates}\cr
   \hbox {in which the volume form $d{\bf v}^\prime$
	     is equal to (1.1))}\,.\cr
	                }
	                   \eqno (1.2a)
		     $$

	               $$
      b)\quad	\Delta_0\sqrt f=
	\sum_{i=1}^n{\p^2 \sqrt f\over \p x^i\p\theta_i}=
	            0\,,\quad 
	\hbox {(BV master-equation for the master-action
	 $S=\log\sqrt f$)}\,.
	                      \eqno (1.2b)
                      $$
		      
		      $$			      
	c)\quad \Delta^2_{d{\bf v^\prime}}=0\,,
	                       \eqno (1.2c)
		      $$ 
obey to the relation
$a)\Rightarrow b)\Rightarrow c)$
and under some assumptions 
they are equivalent (see the details below).

 {\bf Statement 2.}  
    
  {\it The integrand of BV partition function is the semidensity
   $\sqrt{f\dv}$
   which is the natural integration object over Lagrangian
   surfaces in an special odd symplectic superspace.
   In the case if conditions} (1.2) {\it are fulfilled,
    the corresponding integral does not change under
    small variations of Lagrangian surface (gauge-independence
     condition). To the semidensity $\sqrt{f d{\bf v}}$
     corresponds the cohomological class of Lagrangian surfaces.}
     
     The explicit formula for this semidensity was 
     delivered in [\JMPd].        	   	   	
	
     The complete analysis of these statements
    in  the [\SchCMPa],  was particularly founded
     on the relations established in this paper
   between differential forms
    on an usual space and corresponding volume
   forms on an special odd symplectic superspace which
   is associated to its cotangent bundle.
    Focusing our attention  on these considerations
    we come to the construction of $\Delta$ operator, acting
     on semidensities in general odd symplectic superspaces,
     without additional volume form structure.
      This leads us to more clear interpretation
      of the relations in the Statements 1 and 2. In particular
      the condition (1.2b) in terms of semidensities
       receives its invariant formulation and the difference between
        conditions (1.2a,b,c) is formulated exactly.
      But more important is that our considerations
      provide us by the new outlook on
      the invariant semidensity [\CMP,\But]  and
       lead us to construction of
       new integral invariants on embedded surfaces.
      
  We recall shortly the problem of invariant densities
  construction in an odd sympelctic superspace.
  
   In the case if we consider the volume form
   not only on the space (superspace) but
   on arbitrary embedded surfaces we come to 
 the concept of densities on embedded surfaces.
 
  The density of weight $\s$ and rank $k$
  on embedded surfaces is a function
  
  \noindent 
  $A(z,{\p z\over\p\zeta},\dots,{\p^k z\over\p\zeta\dots\p\zeta})$
   which is defined on parametrized surfaces $z(\zeta)$,
   depends on first $k$ derivatives of $z(\zeta)$
   and is multiplied on the $\s$-th power of the determinant
(Berezinian) of surface reparametrization. On the every given
 surface it defines the $\s$-th power of volume form.
  In particular such a concept of density is very useful 
  in supermathematics where the notion of differential forms
  as integration objects is ill-defined [\Gayduk].
  
   In usual mathematics,
  for every $2k$-dimensional surface $C^{2k}$ embedded
   in a symplectic space,
    so called Poincare-Cartan integral invariants
    (invariant volume forms on embedded surfaces)
     are given by the formula
                       $$
   \int_{C^{2k}} \underbrace {w\wedge\dots\wedge w}_{k-{\rm times}}=
              \int 
                 \sqrt
                   {
		  \det 
		  \left
		     (
		{\p x^\mu(\xi)\over\p\xi^i}   
		       w_{\mu\nu}        
	{\p x^\nu(\xi)\over\p\xi^j}
	          \right
		    )
		    }d^{2k}\xi\,,
		                  \eqno (1.3)
		  $$
 where the non-degenerated closed two-form 
 $w=w_{\mu\nu}dx^\mu\wedge dx^\nu$ defines symplectic structure,
  and the functions 
 $x^\mu=x^\mu(\xi^i)$ define some parametrization		  
 of the surface $C$.
		  
 In the case of even symplectic superspace,
  the l.h.s. of (1.3) is ill-defined but  the r.h.s.  of
  this formula  can be straightforwardly generalized,
   by changing determinant on the {\it Berezinian}
    (superdeterminant).
    The properties
    of the integral invariant do not change drastically.
   In particular one can prove  that
     the integrand in the (1.3)
     in the case of even symplectic structure
      in the superspace is
      total derivative 
       and all invariant densities
   on surfaces are exhausted by (1.3)
   as well as in the case of usual symplectic structure 
   [\Poin,\BarScha,\BarSchb]. 
     
    The situation is less trivial in the case of odd symplectic
     space.

     In [\CMP,\But] was analyzed the problem of invariant densities
     existence on the non-de\-ge\-ne\-ra\-ted surfaces embedded in an
      special odd symplectic superspace.
     
         It was proved that there are no invariant densities
   of the rank $k=1$ (except of the volume form itself),
    and it was constructed the semidensity of the rank $k=2$
     which is defined on non-degenerated surfaces of the codimension
      $(1.1)$ (see formulae (3.1, 3.4) below or [\CMP,\But] for details).

 This semidensity takes odd values.
 It is an exotic analogue of Poincar\'e--Cartan
 invariant: the corresponding density
 (the square of this odd semidensity)
 is equal to zero, so it cannot be  integrated
 nontrivially over supersurfaces.
       (The analysis of this semidensity  performed in
      [\CMP] showed that it can be considered as an analog
      of the mean curvature of hypersurfaces in the Riemanian 
    space.)

      Moreover, it was also proved that at least in an normal special
      odd symplectic superspace this odd semidenstiy
      is unique (up to multiplication
       by a constant) in the class of densities of the rank $k=2$
       which are defined on non-degenerated 
       surfaces of arbitrary dimension
        [\But]. 
    
      These results indicate that one have to search
      non-trivial integral invariants 
      (invariant densities  of the weight $\s=1$)
     in higher ($k\geq 2$) order derivatives.
    The tedious calculations which lead to
     the construction of this odd semidensity in the papers [\CMP,\But]
     did not give hope to
    go further for finding them,
     using the technique which was used in these papers.

      The analysis performed in this paper
      shows the important role of semidensity in
      an odd symplectic superspace, revealing
      its meaning in terms
      of differential forms on underlying space.
      It is semidensity in the ambient odd
      symplectic space, not the volume form, 
      which naturally induces invariant densities
      on embedded surfaces (see the Lemma in the Section 3).
       This makes the fact that
       the simplest invariant density on surfaces is the  
       odd semidensity less surprising.
       
       On the other hand our approach allows to construct
       new densities 
       depending on fourth
        order derivatives on surfaces
        embedded in an special odd symplectic superspace.
	
    In the second section we recall the basic definitions
of an odd symplectic and special odd symplectic superspace,
the properties of $\Delta$ operator acting on functions,
 and we give the definition of the
  $\Delta$-operator acting on semidensities
   in an odd symplectic superspace.
  Analyzing these constructions in terms of
   underlying space geometry we come to
   more clear interpretation of BV-formalism geometry.
  
      In the third section   
  we come to natural interpretation
of the odd invariant semidensity on $(1.1)$-codimensional surfaces,
 expressing this semidensity straightforwardly 
 in terms of semidensity in the
  ambient symplectic space. Using this relation and
   operator $\Delta$ on semidensities 
   we come to the main result of this paper
   constructing another semidenstiy
   and  two densities (integral invariants), even and odd,
    of the rank $k=4$
   on  $(1.1)$-codimensional surfaces.
   May be these densities are the simplest (having the lowest rank)
   non-trivial integral invariants on surfaces in an 
  special odd symplectic surfaces.

 \bigskip
 \centerline {\bf 1. $\Ds$ on Semidensities}  
   \medskip
   
         Let $E^{n.n}$  be a superspace with coordinates
        $z^A=$  $x^1,\dots,x^n$, $\theta_1,\dots,\theta_n$;
        $p(x^i)=0,p(\theta_j)=1$, where $p$ is a parity
	\footnote{$^*$}	
	{\noindent We use the definition ,
	 when a point
	  of superspace $E^{n.n}$ is $\Lambda$-point---
  $2n$--plet $(a^1,\dots,$ 
  $a^n,\alpha^1,\dots,\alpha^n)$, where
  $(a^1,\dots,a^n)$ are even and $(\alpha^1,\dots,\alpha^n)$ are
   odd elements of an arbitrary Grassmann algebra $\Lambda$.
  (It is the most general definition
  of superspace suggested by D. Leites and A.S. Schwarz
  as the functor on the category of Grassmann
  algebras [\leites,\SchCMPb].)}.

        We say that this superspace is odd symplectic superspace
       if it is endowed with an odd symplectic structure, i.e.,
       if an odd  closed non-degenerate $2$-form:
       $\Omega=\Omega_{AB}(z)dz^A dz^B$,
       $p(\Omega)=1$, $d\Omega=0$
  is defined  on it [\Berezin, \leitb].

  In the same way as in the standard symplectic calculus one can
 relate to the odd symplectic structure the odd Poisson
      bracket (Buttin bracket) [\Berezin,\leitb,\leites]:
                       $$
                     \{f,g\}=
          {\p f\over\p z^A}
                     (-1)^
                   {fA+A}
                          \Omega^{AB}
              {\p g\over\p z^B}\,,
                       $$
     where $\O^{AB}=\{z^A,z^B\}$
 is the inverse matrix to $\O_{AB}$ :
$\O^{AC}{ \O_{CB}}=\d^A_B$ ,

\noindent($\O^{AB}=-\O^{BA}(-1)^{(A+1)(B+1)}$).

 To a function $f$ there corresponds
  the Hamiltonian vector field
                       $$
         {\bf D}_f=\{f,z^A\}
      {\p\over\p z^A}\,,
                 \quad
         {\bf D}_f(g)=\{f,g\},\quad
          \O({\bf D}_f,{\bf D}_g)=-\{f,g\}\,.
                                     \eqno (2.1)
                      $$
    (See for the details e.g. [\CMP].)
		      
    The condition of the closedness of the form defining
     symplectic structure leads to the Jacoby
    identities:
                       $$
           \{f,\{g,h\} \}(-1)^
              {(f+1)(h+1)}+
                   {\rm cycl.\quad permutations}=0
		                  \eqno (2.2)
                        $$
    Using the analog of Darboux Theorem [\Shander]
   one can consider the coordinates in
    which the symplectic structure  and the
   corresponding Buttin bracket
    have locally the canonical expressions:
    $\O=\sum dx^i d\theta_i$
    and respectively                     
                        $$
         \{x^i,x^j\}=0,\, \{\theta_i,\theta_j\}= 0,\,
                 \{x^i,\theta_j\}=\d^i_j\,,
                      \{f,g\}=
                      \sum_{i=1}^n
                       \left(
         {\p f\over\p x^i}
         {\p g\over\p \theta_i}
                           +(-1)^f
         {\p f\over\p \theta_i}
         {\p g\over\p x^i}
                       \right)\,.
                                           \eqno (2.3)
                         $$
   We call these coordinates {\it Darboux coordinates}.

   The odd symplectic space provided with a volume form
                         $$
	 d{\bf v}=\rho(z)dz^1\dots dz^{2n}
	                                 \eqno (2.4)
			$$
will be called {\it special odd symplectic superspace}.
We suppose that the volume form is non-degenerated, i.e.
for the every point $z_0$ the number (non-nilpotent) part of $\rho(z_0)$
is not equal to zero.

In the case if in the special odd symplectic superspace
there exist Darboux coordinates such that the volume form
 in these coordinates is given by (1.1),
 then this space will be called
{ \it normal special odd symplectic superspace}.

The action of $\Delta$ operator on an arbitrary
functions in the special odd symplectic
superspace is equal (up to coefficient)
to the divergence w.r.t. volume form (2.4) of the Hamiltonian
vector field corresponding to this function [\JMPa,\MPL].
Using (2.1) we come to  the formula
                              $$
	\D f={1\over 2}(-1)^{f}{\rm div}_{d{\bf v}}{\bf D}_{f}=
	         {1\over 2}(-1)^{f}
		 \left(
		 (-1)^{{\bf D}_{f}A+A}
		   {\p\over\p z^A}
		     \{f,z^A\}+
		       D_f^A{\p \log\rho(z)\over\p z^A}
                	\right)\,.
			           \eqno (2.5)	            
	                      $$
			      		      
 In Darboux coordinates:
                        $$
	\D f=\Delta_0 f+ {1\over 2}\{\log\rho, f\}\,,		           
	               \eqno (2.6)
			$$	     
 where $\rho(z)$ is given by (2.4), and
                   $$
     \Delta_0 f=\sum_{i=1}^n{\p^2 f\over \p x^i\p\theta_i}\,.
                         \eqno (2.7)
                   $$			 		   			
(Later on we use mostly only Darboux coordinates.)
				
This operator obeys to the relations [\BVb,\MPL]:			
			$$
	\D \{f,g\}=\{\D f,g\}+(-1)^{f+1}\{f,\D g\}\,,
                        $$
   			$$
	\D (f\cdot g)=\D f\cdot g+(-1)^f f\cdot\D g+(-1)^{f}\{f,g\}\,.
	                            \eqno (2.8)
                        $$

\smallskip

Now we define the action of operator $\Ds$ on semidensties
in an odd symplectic superspace.

{\bf Definition} 
{\it If ${\bf s}$ is a semidensity in an odd symplectic superspace
 and $s(z)|dz|^{1/2}$ is local expression
for this semidensity in Darboux coordinates 			
$z^A=(x^1,\dots, x^n,\theta_1,\dots,\theta_n)$ then the local
 expression for the density $\Ds {\bf s}$
 in these coordinates is given by the following formula:}
                     $$
              \Ds {\bf s}=(\Delta_0 s(z))|dz|^{1/2}=
	 \sum_{i=1}^n{\p^2 s\over \p x^i\p\theta_i}
	          |dx d\theta|^{1/2}\,.
	                              \eqno (2.9)
	             $$	
		     		                
One can prove that the r.h.s. of this formula defines
 the density also, i.e. if 
 $\t z=(\t x^1,\dots,\t x^n,\t\theta_1,\dots,\t\theta_n)$ are another
 Darboux coordinates then
                   $$
		   \left(
  \sum_{i=1}^n{\p^2 \over \p x^i\p\theta_i}
                 s(z)
		 \right)_{z(\t z)}
                  \cdot
		    \left(                   
	 {\rm Ber}{\p z(\t z)\over\p\t z}
	        \right)^{1/2}
                    =
               \sum_{i=1}^n
      {\p^2\over \p \t x^i\p\t\theta_i}
	          \left(
                  s(z(\t z))
                    \cdot
      {\rm Ber}{\p z(\t z)\over\p \t z}^{1/2}
                      \right)\,.		      
		  		\eqno (2.10) 
		   $$	
  Canonical transformations 
  from some Darboux coordinates
  to another Darboux coordinates infinitezimally
  are generated by an odd
  function (Hamiltonian) via corresponding Hamiltonian vector filed (2.1).
  To an odd function $Q(z)$ corresponds transformation
   $\t z^a=z^A+\varepsilon \{Q,z^A\}$.
 To the action of this transformation on the semidenity ${\bf s}$
     corresponds differential
    $\d_Qs=\Delta_0 Q\cdot s-\{Q,s\}$,
    because $\d s=-\varepsilon\{Q,s\}$ and 
    $\d |dz|=\varepsilon\d {\rm Ber}(\p z/\p\t z)|dz|=
    2\Delta_0Q|dz|$.
    Using that $\Delta_0^2=0$ and (2.8) we come to commutation relations
   $\Delta_0\d_Q=\d_Q\Delta_0$. This leads to relation (2.10),
   which proves the correctness of definition (2.9). 
   One can check the relations (2.10) by straightforward computations also,
   using the properties of the operator $\Delta_0$ 
   which were investigated in details in [\BVb].

The action of differential $\d_Q$
on semidensities can be rewritten in a explicitly invariant way:
                  $$
		  \d_Q{\bf s}=Q\cdot\Ds{\bf s}+\Ds( Q{\bf s})=
		        [Q,\Ds]_+{\bf s}\,.
			              \eqno (2.11)
		 $$		      

Contrary to the operator $\D$ on functions, the operator $\Ds$
on semidensities does not need volume structure.

On a special odd symplectic superspace
we can construct new invariant objects,
expressing them via the semidensity related with volume form
(${\bf s}=\sqrt\dv$) and operator $\Ds$:
                       $$
      {\bf s}=\sqrt\dv \,
      \hbox {---semidensity ($\sigma={1\over 2}$)}\,,
                                   \eqno (2.12a)
			$$
			$$	   
		\Ds {\bf s}=\Ds\sqrt\dv\,
		\hbox {---semidensity ($\sigma={1\over 2}$)}\,,
		                \eqno (2.12b)
			$$
			 $$	
	{\bf s}\Ds{\bf s}=\sqrt\dv\Ds\sqrt\dv\,
	\hbox {---density ($\sigma=1$)}\,,
	                       \eqno (2.12c)
			$$
			$$       
	         {1\over{\bf s}}\Ds{\bf s}= 
		 {1\over\sqrt\dv}\Ds\sqrt\dv\,	        			 
		 \hbox {---function ($\sigma=0$)}\,.
		                          \eqno (2.12d)
			 $$

Using (2.6)---(2.9) we can see that operator
$\Ds$ obeys to the following properties:
                      $$
		      (\Ds)^2=0\,,
		       $$
                      $$
		    \Ds (f\cdot \sqrt\dv)=
		    (\D f)\cdot \sqrt\dv+
		    (-1)^{f}f\cdot \Ds \sqrt\dv\,,
		                        \eqno (2.13) 
                    $$      
  and
                 $$
		 \D^2 f=
		 \{{1\over\sqrt\dv}\Ds\sqrt\dv,f\}\,.
		                   \eqno (2.14)
			$$
 To clarify the geometrical meaning of the definition
 (2.9) and of its correctness
 we consider the following example

 {\bf Example 1} Let $M$ be a $n$-dimensional space and
 $T^*M$ be its cotangent bundle.
 Let $ST^*M$ be a superspace associated with
 $T^*M$. To local coordinates $(x^1,\dots,x^n)$
 on $M$ correspond the local coordinates 
 $z^A=(x^1,\dots, x^n,\theta_1,\dots,\theta_n)$  
  on $ST^*M$. 
  
  The odd coordinates $\theta_j$ transforms
  via the differential of corresponding
   transformation of the even coordinates $x^i$:
   
                $$
\t x^i=\t x^i(x),\quad\t\theta_i=\sum_{k=1}^n
      {\p x^k(\t x)\over\p \t x^i}\theta_k\,.	
		  		     \eqno (2.15)
	        $$
Therefore one can define the canonical odd symplectic structure
on $ST^*M$ in a such way that 
$z^A=(x^1,\dots, x^n,\theta_1,\dots,\theta_n)$				     
are Darboux coordinates for this symplectic structure.
The pasting formulae (2.15) provide the correctness of the
definition of this structure. 

The relations between the cotangent bundle structure on $T^*M$
and the odd canonical symplectic structure (2.3) on $ST^*M$
are based on the canonical map $\tau_{_M}$:
               $$
	       \tau_{_M}
	       \left(
	       T^{i_1\dots i_k}
	       {\p\over\p x^{i_1}}
	       \wedge\dots\wedge
	       {\p\over\p x^{i_k}}
	           \right)=
	       T^{i_1\dots i_k}
	       \theta_{i_1}\dots\theta_{i_k}\,,
	              $$
 between polyvectorial antisymmetric fields on $M$ and
 functions on $ST^*M$. This map transforms the
 Schoutten bracket to the odd canonical Poisson (Buttin)
 bracket (2.3) [\leites,\buttin]:
                        $$
			\tau_{_M}
			\left(
			\left[
			{\bf T_1},{\bf T_2}
			\right]
			  \right)=
			\{
		\tau_{_M} \left({\bf T_1}\right),
		\tau_{_M} \left({\bf T_2}\right)
		          \}\,.
			           \eqno (2.16)
			 $$
 
 Now we consider the  map $\tau_{_M}^{\#}$
  from the differential forms on $M$
  to the semidensities on $ST^*M$,
  defining it in the following way:
		
                   $$
		 \eqalign
		 {
	\tau_{_M}^{\#}(1)&=\theta_1\dots\theta_n|dx d\theta|^{1/2},\cr
	\tau_{_M}^{\#}(dx^i)&=(-1)^{i+1}
 \theta_1\dots
 \widehat\theta_i\dots\theta_n|dx d\theta|^{1/2}\,,\cr	 
\tau_{_M}^{\#}(dx^i\wedge dx^j)&=(-1)^{i+j }
    \theta_1\dots
    \widehat\theta_i\dots\widehat\theta_j
    \dots\theta_n|dx d\theta|^{1/2}\,,\,(i<j),\cr	 
        &\dots\cr
  \tau_{_M}^{\#}(dx^{i_1}\wedge\dots\wedge dx^{i_k})&=
           (-1)^{i_1+\dots+i_k+k}
      \theta_1\dots
    \widehat\theta_{i_1}\dots\widehat\theta_{i_k}
    \dots\theta_n|dx d\theta|^{1/2},\,
     (i_1<\dots <i_k), \cr	 
                 }
		$$		
                 $$
     \tau_{_M}^{\#}(f(x)w)=f(x)\tau_{_M}^{\#}(w)\,,\quad
  \hbox{for every function $f(x)$ on $M$}\,,
                            \eqno (2.17)
                               $$
where the sign $\,\,\widehat {}\,\,$
 means the omitting of corresponding term.			     
	For example
if $M$ is two-dimensional space, then
$\tau^{\#}_{_M}(f(x))=f(x)\theta_1\theta_2$,
$\tau^{\#}_{_M}(w_1(x)dx^1+w_2(x)dx^2)=
w_1\theta_2-w_2(x)\theta_1$,
$\tau^{\#}_{_M}(w(x)dx^1\wedge dx^2)=-w(x)$
\footnote
{$^*$}{The square of this map 
$\tau^{\#}\colon\, w\rightarrow (\tau^{\#}(w))^2$
transforms differential forms on $M$ to {\it density}
(volume form) on $ST^*M$. 
This map was constructed in [\SchCMPa]
via the superspace $STM$ associated to tangent bundle
$TM$ and additional arbitrary volume form on $ST^*M$.}.

We say that semidensity ${\bf s}$ corresponds to
 differential form $w$ (to the linear combination of
  differential forms $\sum w_k$) if ${\bf s}=\tau^{\#}_M (w)$
  (${\bf s}=\tau^{\#}_M (\sum w_n)$).

The correctness of (2.17) follows from the
 fact that for transformations (2.15)
              $$
	    \det {\p x\over\p\t x}=
	           {\rm Ber}
	           \left(
		{\p(x,\theta)\over \p(\t x,\t\theta)}
		\right)^{1/2}  \,,          
		    \qquad
		    \left(
                    {\rm  Ber}
                    \pmatrix
                 {A&B\cr C&D\cr}
                        =
           {\det\, (A-BD^{-1}C)\over \det\, D}
	               \right)\,.
		        \eqno (2.18)
                        $$
    
 The action of operator $\Ds$ corresponds to
 the action of exterior differential:
                            $$ 		  
	\Ds\circ\tau_{_M}^{\#}=\tau_{_M}^{\#}\circ d\,,
	                \qquad {\rm and}\qquad
     \tau_{_M}^{\#}({\bf T}{\cal c} w)=
    \tau_{_M}({\bf T})\cdot\tau_{_M}^{\#}(w)\,,
             \eqno (2.19)
	        $$        	  
where ${\bf T}{\cal c} w$ is the inner product of
polyvectorial field ${\bf T}$ with differential form $w$.

 If the semidenstity
 ${\bf s}$ in $ST^*M$ corresponds to the volume form
 (differential $n$-form) $w$ 
 on $M$ and the structure
 of special odd sympelctic superspace on $ST^*M$ 
 is defined by the square of this semidenstiy
 ($d{\bf v}={\bf s}^2$)
 then  the action of $\D$ corresponds to the
 divergence w.r.t. to the volume form $w$
 on $M$: 
 $\D\circ\tau_{_M}=
 \tau_{_M}\circ {\rm div}_w$.
 (See also [\JMPd,\SchCMPa].)

 Of course the essential difference of the odd 
 symplectic superspace $ST^*M$ from
 the cotangent bundle $T^*M$ is that 
  canonical transformations
 in $ST^*M$ are not exhausted by (2.15).
 For example if 
  we consider the action of Hamiltonian
  $Q=L\theta_1\dots\theta_n$ 
  on the differential form $w=dx^1\wedge\dots\wedge dx^n$
  then we obtain using (2.11) that
  $\d w=(\tau_{_M}^{\#})^{-1}\d_Q\tau_{_M}^{\#}w=dL$.
   Only in the special case if an odd Hamiltonian $Q$ corresponds
 to the vector field $T^i{\p\over\p x^i}$
 ($Q=T^i\theta_i$), then this Hamiltonian
 induces infinitezimal canonical
 transformation, which corresponds to the point
 transformation on $M$:
 $(\tau_M^{\#})^{-1}\d_Q\tau_M^{\#} w$ is equal to Lie
 derivative of $w$ along the vector field 
 $T^i{\p\over\p x^i}$.
  
 The initial space $M$ is the Lagrangian 
  $(n.0)$-dimensional surface in $ST^*M$.
 If $L^{n.0}$ is an arbitrary $(n.0)$-dimensional 
  Lagrangian surface in $ST^*M$ then
   there exist  Darboux coordinates 
   such that $L^{n.0}$ is given in 
   these coordinates locally by the
    equations $\theta_1=\dots=\theta_n=0$.	         	         	         
A general canonical transformation transforms the initial
 space $M$ to some $(n.0)$-dimensional
 Lagrangian surface in $ST^*M$.
 An arbitrary semidensity ${\bf s}$
  on $ST^*M$
 corresponds by (2.17) to the linear combination of
 differential forms on $L^{n.0}$, 
 via the map $\tau^{\#}_{_L}$.
 (In the case if we consider the points of superspace as
  $\L$-points  where $\L$ is an arbitrary Grassmann algebra,
  (see the footnote in the beginning of the section)
   then one have to consider differential forms with coefficients
   in this algebra $\L$.)
  
  \smallskip
   
  This example is basic one, because
   locally every odd symplectic superspace 
   can be considered as superspace $ST^*L$
   associated to the cotangent bundle $T^*L$
   of its $(n.0)$-dimensional Lagrangian surface $L$.
   The semidensity can be integrated
   over arbitrary $(n-k.k)$-dimensional
   Lagrangian surface in an odd symplectic superspace
    (See [\SchCMPa] and [\JMPd] for explicit formula).
    In the case if the Lagrangian surface is $(n.0)$-dimensional,
  the integral of semidensity ${\bf s}$ 
  over this surface $L$
   is nothing but the integral of corresponding
differential form $(\tau_{_L}^{\#})^{-1}{\bf s}$
 by this surface.
 
   Now using these constructions
   we return to statements (1.2) of BV formalism
   geometry. 
   
    One have straightforwardly deal with an odd symplectic 
    superspace provided with semidensity
     ${\bf s}=\sqrt\dv=(f|dx d\theta|^{1/2})$, which
     corresponds to the exponent of master-action
     $f=e^{2S}$.

   The master-equation (1.2b) in terms of the operator
    $\Ds$ can be rewritten as the condition of
    the semidensity (2.12b) vanishing: $\Ds{\bf s}=0$. 
     It means according to (2.19),
      that the differential
     form $(\tau^{\#})^{-1}{\bf s}$ corresponding 
     to the semidensity $\sqrt{d{\bf v}}$
     is closed. 
     
     Under the infinitezimal transformation 
     to another Darboux coordinates,
 according (2.11) this semidenity transforms as
                        $$
	\d_Q {\bf s}=\Ds(Q{\bf s})\,,
	\eqno (2.20)
		                $$
 and corresponding form 
 $w=({\tau^{\#}})^{-1}{\bf s}$ changes on the exact
 form $d({\tau^{\#}}^{-1}(Q{\bf s}))$.
 
    If in some Darboux coordinates this semidensity is expressed by 
    the formula
            $$
 {\bf s}=s(x,\theta)|dxd\theta|^{1/2}=
     (a(x)+a^i(x)\theta_i+\dots+ c\theta_1\dots\theta_n)
     |dxd\theta|^{1/2}
                          \eqno (2.21)
	  $$
then $c$ is a constant, because
 the form $({\tau^{\#}})^{-1}{\bf s}$ is closed.
 From (2.20) it follows that
  this constant does not depend on the choice 
  of Darboux coordinates. In fact it characterizes cohomological
  class of the form $w=(\tau^{\#})^{-1}{\bf s}$  		  	                
  on the $(0.n)$-dimensional Lagrangian plane
  (see [\SchCMPa]):
                        $$
		c= \int_{L^{0.n}}{\bf s}=\int_
              {x^1=x^1_0,x^2=x^2_0,\dots,x^n=x^n_0}
	           s(x,\theta)d\theta_1\dots d\theta_n\,.
		                 \eqno (2.22)    
			$$
				 
     One can say more: the condition of closedness of the form
     corresponding to the semidensity ${\bf s}$ means
     that by choosing appropriate Hamiltonian $Q(z,t)$,
     and integrating the relation (2.20) one comes to
     the canonical transformation
      to new Darboux coordinates in which:
                         $$
	{\bf s}=s(x,\theta)|dxd\theta|^{1/2}=
     \left(1+ c\theta_1\dots\theta_n\right)
     |dxd\theta|^{1/2}\,.
                          \eqno (2.23)
			$$
 We find this canonical transformation
 in the way similar to [\SchCMPa]
  using the correspondence between semidensities
 and differential forms and 
 the Principal Formula of Differential Forms Differential Calculus
 [\Gui].
 
 Let ${\bf s}_0$ and ${\bf b}$
 be arbitrary "closed" semidensities,
 (i.e.,
  the corresponding differential forms are closed).
       We consider the odd function $Q$
 which obeys to relation 
 $Q{\bf s}_t=-{\bf b}$, where 
 ${\bf s}_t={\bf s}_0+t\Ds{\bf b}$.
 (We suppose that the semidensity ${\bf s}_t$ is not degenerated,
  i.e.  $n$-form
   corresponding to ${\bf s}_t$ is not equal to zero).
 One can see that canonical transformation $z=z(\t z, t)$
 which is defined by the equation:
                        $$
		\cases
		{
	 {dz(\t z,t)\over dt}=\{Q,z(\t z,t)\}\cr
	   z(\t z,t)\vert_{t=0}=z\cr
			  }\,,
			  $$
transforms the semidenstiy ${\bf s}_t$ 
to the semidensity ${\bf s}_0$
for an arbitrary $t$.
 This follows from initial conditions and from the fact that
                    $$
           {d {\bf s}_t\over dt}
           =\Ds{\bf b}+\d_Q {s}_t=
	   \Ds({\bf b}+Q{\bf s}_t)=
                     0\,.
		     $$
We put
${\bf s}_0$ to be equal to the semidenstiy (2.23) and choose 
the semidensity ${\bf b}$
such that ${\bf s}_1$ is equal to the semidenstiy (2.21).
(This is possible, because locally every closed $n$-form is exact,
if $n\geq 1$.)			  
			   
The condition of $c\not=0$ in (2.21--- 2.23)
is the obstacle to the condition (1.2a).

 Now we analyze the condition (1.2c). 
  From (2.14) it follows that the condition
     (1.2c) means that the function (2.12d)
     is equal to an odd constant $\nu$, and
     $\Ds {\bf s}=\nu{\bf s}$. 
     One can see using correspondence between semidensities
     and differential forms that
     all the  solution to this equation
    are  ${\bf s}=\Ds {\bf h}-\nu {\bf h}$,
   where ${\bf h}$ is an arbitrary semidensity.
The odd constant $\nu\not=0$
     is the obstacle to the condition (1.2b).
   
 We come to the
   
   {\bf Proposition}
   {\it In the case
   if the odd symplectic superspace is provided by the volume form
   $\dv$, such that $\D^2=0$,
   then to the volume form $\dv$ corresponds the odd constant
   $\nu$: $\Ds\sqrt\dv=\nu\dv$ and the semidensity 
   $\sqrt\dv$ has the form
   $\Ds {\bf h}-\nu {\bf h}$.
   If the odd constant $\nu$ is equal to zero, then
   the master-equation $\Ds\sqrt\dv=0$ holds. 
   In this case to the volume form $\dv$ corresponds the constant
   $c$ which is equal to the integral of semidensity
   $\sqrt\dv$ over $(0.n)$-dimensional Lagrangian plane} (2.22).
   {\it The volume form in this case can be reduced 
   to the form} (2.23).
   {\it In the case if this constant $c=0$ then
   the superspace is the normal special odd symplectic superspace.}
       
   \smallskip

    This Proposition removes uncorrectness of the
     considerations about equivalence of conditions
     (1.2) which was done in [\JMPd] and [\SchCMPa].
     On the other hand the statements of this Proposition
     in non explicit way was contained in
     the statements of Lemma 4 and Theorem [5] of the paper
      [\SchCMPa].
     The analysis performed here in terms of $\Ds$
     operator is essentially founded on these results.

\bigskip

  \centerline{\bf 3. Invariant semidensity on $(1.1)$-
  codimensional surfaces}
  
  \medskip
  
  In this section first we 
   recall explicit formulae for
    the odd invariant semidensity on
    non-degenerated surfaces of codimension $(1.1)$
    embedded in an special odd symplectic superspace
     ([\CMP,\But]). Then we rewrite this semidensity
      in terms
     of the semidensity $\sqrt\dv$ of the 
     ambient superspace and suggest
      the construction of pull-back of arbitrary
      semidensity from the ambient odd symplectic superspace
       on embedded $(1.1)$-codimensional surfaces.
       Using the semidensity
     $\Ds\sqrt\dv$ we will construct the new densities
     on embedded non-degenerated surfaces.
     
     The surface is called  non-degenerated if the
     sympelctic structure of the space generates
     the symplectic structure on the embedded surface also,
     i.e. if the pull-back of the symplectic $2$-form 
     on the surface is non-degenerated
      $2$-form. As usual we call this symplectic structure 
      on an 
      embedded surface {\it induced symplectic structure}. 
     
          Let $z^A$ be Darboux coordinates
	on an special odd symplectic superspace
	$E^{n.n}$
	with volume form $\dv=\rho(z)|dz|$.
	It is convenient here to use for Darboux coordinates
	notations
	$z^A=(x^\mu,\theta_\mu),$ $(\mu=(0,i)=$ $(0,1,\dots,n-1)$,
	$i=(1,.\dots,n-1))$.
       Let $z(\zeta)$ be an arbitrary parametrization
       of an arbitrary non-degenerated 
       surface of codimension $(1.1)$, embedded
        in this  special odd symplectic superspace.
	($\zeta=(\xi^i,\eta_j)$, $\xi^i$ and $\eta_j$
	are even and odd parameters respectively, 
	($i,j=1,\dots,n-1$).

  As it was mentioned in Introduction there is no non-trivial
  invariant density of the rank $k=1$
   on non-degenerated surfaces.
   
      The invariant semidensity of the rank $k=2$
      (depending on first and second derivatives of $z(\zeta)$)
       is given by the 
      following formula:
                    $$
              A\left(z(\zeta),{\p z\over\p\zeta},
    {\p^2 z\over\p\zeta\p\zeta}\right)|d\zeta|^{1/2}=
                    $$
		    $$    
	           \left(
	\Psi^A{\p \log\rho(z)\over\p z^A}\,
  -\,\Psi^A\O_{AB}{\p^2 z^B\over\p\zeta^\a\p\zeta^\beta}
            \O^{\a\beta}(z(\zeta))
                  (-1)^
       {B(\a+\beta)+\a}
               \right)
	       |d\zeta|^{1/2}\,,
                                    \eqno (3.1)
                   $$
  where $\O_{AB}dz^A dz^B$ is the two-form defining the odd
  sympelctic 
   structure on $E^{n.n}$
   and $\O^{\a\beta}$ is the tensor inversed to 
   the two-form which defines induced symplectic structure on the
   surface.
   The vector field ${\bf \Psi}$ is defined as follows:
   one have to consider 
   the pair of vectors $(\bf H,\bf\Psi)$,
   $\bf H$-even and $\bf\Psi$-odd, 
   which are symplectoorthogonal to the surface
    and obey to the following conditions:
                    $$
      \O\left({\bf H,\bf\Psi}\right)=1,\,\,
        \O\left({\bf\Psi,\Psi}\right)=0
	 \quad ({\rm symplectoorthonormality\,
	 conditions})\,,
	                 \eqno (3.2)
		   $$
		   $$        	 
       d{\bf v}\left(\left\{{\p z\over\p\zeta}\right\},
       {\bf H},{\bf\Psi}\right)=1\quad
        \hbox{(volume form normalization conditions)}\,.
                           \eqno (3.3)
       		    $$
These conditions fix uniquely the 
vector field ${\bf \Psi}$. (See for  details [\CMP]).
           
 The explicit expression for this semidensity
  was calculated in [\But] in  terms
    of dual densities: If $(1.1)$-codimensional
   supersurface $M$ is given not by parametrization, but
    by the equations
      $f=0,\varphi=0$,
   where $f$ is an even function and $\varphi$ is an odd function
 then to the semidensity (3.1) there corresponds the 
 dual semidensity:
                   $$
                 {\tilde A}\Big\vert_{f=\varphi=0}=
           {1\over\sqrt{\{f,\varphi\}}}
                  \left(
                  \D f-
          {\{f,f\}\over 2\{f,\varphi\}}
                \D\varphi-
              {\{f,\{f,\varphi\}\}\over\{f,\varphi\}}-
                  {\{f,f\}\over 2\{f,\varphi\}^2}
              \{\varphi,\{f,\varphi\}\}
                       \right)\,.
                                       \eqno(3.4)
                        $$
  This function 
 is multiplied by the square root of the corresponding
  Berezinian (superdeterminant) under the
   transformation
     $f\rightarrow$ $a f+\alpha\varphi$,
     $\varphi\rightarrow$ $\beta f+b\varphi$,
   which does not change the surface $M$.

 Now we rewrite the semidensity (3.1)
 straightforwardly via the semidensity $\sqrt\dv$
 on the ambient special odd symplectic superspace $E^{n.n}$.
 For a given non-degenerated surface $M$ of codimension
 $(1.1)$ embedded in $E^{n.n}$ one can
 choose Darboux coordinates (in a vicinity of arbitrary point)
 such that in these Darboux coordinates
 the surface $M$ is given by equations
                       $$
		      x^0=\theta_0=0\,.
		       \eqno (3.5)
		      $$
We call these Darboux coordinates adjusted to the surface $M$.		      
We note that if $(x^0,x^i,\theta_0,\theta_j)$
are Darboux coordinates in $E^{n.n}$
 adjusted to the surface $M$, then
 $(x^i,\theta_j)$ are Darboux coordinates
 on the surface $M^{n-1.n-1}$ w.r.t. induced symplectic structure.
We consider the following
parametrization of the surface $M$ in adjusted 
Darboux coordinates:       
                        $$
		       \cases
		          {
			x^0=0,\theta_0=0,\cr
			x^i=\xi^i,\theta_i=\eta_i\,,\,\hbox
			{for $1\leq i\leq n-1$}
			  }
			                  \eqno (3.6)
			 $$		  	
 The conditions 
of symplectoorthonormality in (3.2) give
that  ${\bf H}=a{\p\over \p x^0}+\beta{\p\over\p\theta_0}$
and 
${\bf\Psi}={1\over a}{\p\over\p \theta_0}$,
where $a$ and $\beta$ are arbitrary numbers
 ($a$ even and $\beta$-odd).
The condition of the volume form 
 $d{\bf v}=\rho|dz|$ normalization:
              $$
	      \rho |dx d\theta|\left(
	      {\p\over\p x^1},\dots
	      {\p\over\p x^{n-1}},{\bf H};
	      {\p\over\p \theta_1},\dots
	      {\p\over\p \theta_{n-1}},{\bf \Psi}
	           \right)=
		   \rho a^2
		   =1\,,
		   $$
gives that
                   $$
		   a={1\over\sqrt\rho}\,\quad{\rm and}\quad
  {\bf \Psi}=\sqrt\rho{\p\over\p\theta_0}\,.
		                   \eqno (3.7)
		  $$

Hence (3.1) in these coordinates
 is reduced 
 on $M$ to the semidensity
                    $$
	A=\sqrt\rho{\p\log\rho\over\p\theta_0}
	|\widetilde{dx d\theta}|^{1/2}
	\Big\vert_{x^0=\theta^0=0} =
	    2{\p\sqrt\rho\over\p\theta_0}
	    |\widetilde{dx d\theta}|^{1/2}
	    \Big\vert_{x^0=\theta^0=0}=
	                 $$
			 $$
	        2{\p\sqrt\rho\over\p\theta_0}
	          {\rm Ber}
		    \left(
		  {\p(x^i,\theta_j)\over\p(\xi^i,\eta_j)}
		     \right)^{1/2}
		     |d\xi d\eta|^{1/2}
		     \Big\vert_{x^0=\theta^0=0}\,,
	                    \eqno (3.8)
	             $$
where 
 $\widetilde{dx d\theta}=
 dx^1\dots dx^{n-1} d\theta_1\dots d\theta_{n-1}$.

 The infinitezimal
  transformation from adjusted Darboux coordinates to another
  adjusted Darboux coordinates 
  is generated by an odd Hamiltonian $Q$,
  which according to (3.5)
   subject to the folowing conditions:
                    $$
	 {\p Q(x,\theta)\over\p x^0}
	 \Big\vert_
	{x^0=\theta_0=0}=
	 {\p Q(x,\theta)\over\p \theta^0}
	      \Big\vert_
	{x^0=\theta_0=0}=0\quad\,.
	                           \eqno (3.9)
		 $$

The differential $\d_Q$ (2.11) corresponding
to the action of infinitezimal transformation
generating by the Hamiltonian $Q$ subject to condition (3.9)
on the semidenstiy
 $\sqrt\rho|dz|^{1/2}$
has to commute 
 with the derivative along $\theta_0$  on the surface:
                  $$
 	      \left(
	{\p\over\p\theta_0}\d_Q-
	\widetilde{\d_Q}{\p\over\p\theta_0}
	      \right)
	    \Big\vert_
	{x^0=\theta_0=0}=0\,.
	                  \eqno (3.10)
	   $$
Here $\d_Q$ acts on a semidensity 
in the symplectic space $E^{n,n}$ and 
$\widetilde{\d_Q}$ acts on a semidensity on the surface $M$
provided with induced symplectic structure.
One can check the condition (3.10),
using (2.11) and 
noting that $(x^i,\theta_i)$ are Darboux coordinates 
on the surface $M$ w.r.t.
induced symplectic structure on $M^{n-1.n-1}$.

The
  formula (3.8) defines correctly semidensity on the surface
  $M$ not only for semidenstiy, related with 
 the volume form, but for arbitrary semidenstiy,
 even in the case if it is an odd semidensity, 
 and the corresponding volume form is equal to zero.
  We come to the following statement
 
 {\bf Lemma}
   {\it To every semidensity ${\bf s}$ in the 
   odd symplectic superspace $E^{n.n}$ 
   corresponds semidensity $\A({\bf s})$ 
   defined on non-degenerated
   $(1.1)$-codimensional surfaces embedded in this superspace. 
   
   The value of this semidensity $\A({\bf s})$
   on every $(1.1)$-codimensional
   surface 
   $M$ in Darboux coordinates adjusted to the surface $M$
    is given by the equation}:
                          $$
   \A(M,\,{\bf s})={\p s(x^\mu,\theta_\nu)\over\p\theta_0}
                  \Big\vert_{x^0=\theta_0=0}
	    |\widetilde{dx d\theta}|^{1/2}\,, \quad
	      ({\bf s}=s(x^\mu,\theta_\nu)|dx d\theta|^{1/2})\,.
	                     \eqno (3.11)
			     $$

The formula (3.11) has clear meaning in terms of differential
forms. Consider the Lagrangian surface $L$ in $E^{n.n}$,
which is given in Darboux coordinates $(x^\mu,\theta_\nu)$
 adjusted to the surface $M$
 by the equations $\theta_0=\dots=\theta_{n-1}=0$. The intersection
  $\widetilde L=L\cap M$ of this Lagrangian surface with $M$
  will be Lagrangian subsurface in the symplectic manifold $M$,
   which is given 
    by equations $\theta_1=\dots=\theta_{n-1}=0$ on $M$.
   It is easy to see from (2.17) and (3.11) that if
  the semidensity ${\bf s}$ corresponds to the linear combination
  $\sum w_k$ of differential forms on $L$, 
  ${\bf s}=\tau_{_L}^{\#}(\sum w_k)$, then the semidensity
  (3.11) on $M$ corresponds (up to the sign)
  to the pull-back on $\widetilde L$
  of this differential form:
                       $$
\A(M,{\bf s})=-\widetilde{\tau_{_{\widetilde L}}^{\#}}\circ\iota^*
\circ (\tau_{_L}^{\#})^{-1}(\bf s)\,.
                         \eqno (3.12)
			 $$ 
 $\iota\colon\,\, \widetilde L\hookrightarrow L$ is the embedding map
 and $\widetilde{\tau_{_{\widetilde L}}^{\#}}$ is the
 map (2.17) from differential forms on Lagrangian surface 
 in $M$ to semidensities on $M$.
  			 
\smallskip

The statement of this
Lemma gives us not only the alternative way to prove
the correctness	of the semidensity (3.1),
but allows us to construct semidensity on embedded surfaces via 
odd semidensities on the ambient superspace,
which cannot be yielded from volume forms.

 The odd semidenstiy (3.1) is nothing but
 $\A(M,\,\sqrt\dv)$. 
 To the semidensity $\Ds\sqrt\dv$
 in the special odd symplectic superspace with volume form $\dv$ 
 there corresponds the even
 semidensity $\A(M,\,\Ds\sqrt\dv)$ which cannot be represented
 (3.1)--like, because the square of the odd semidensity 
 $\Ds\sqrt\dv$ is equal to zero. 

 These semidensities
 correspond to differential forms on 
  Lagrangian $(n-1.0)$-di\-men\-si\-onal
   subsurfaces in $M$, hence they can be
  integrated  over them.
   (Moreover they can be integrated 
  over arbitrary $(n-1-k.k)$-dimensional
  Lagrangian subsurfaces in $M$ also (see [\SchCMPa,\JMPd]).)

 One can construct 
 densities of weight $\s=1$
 via the semidensities  $\A(M,\,\sqrt\dv)$ and
 $\A(M,\,\Ds\sqrt\dv)$: 
 
                  $$
		  P_0=\A^2(\Ds\sqrt\dv)\,\quad{\rm and}\quad
		  P_1=\A(\sqrt\dv)\A(\Ds\sqrt\dv)\,.
		       \eqno (3.13)
		     $$
The density  $P_0$ takes even values, 
the density $P_1$ takes odd values.
In general case these
densities give non-trivial integration objects
over non-degenerated $(1.1)$-codimensional surfaces
embedded in an special odd symplectic superspace with
volume form $\dv$.

The densities $P_0$ and $P_1$ have rank $k=4$
(i.e. depend in general 
on derivatives of the parametrization $z(\zeta)$
up fourth order).
It follows from the fact that 
the semidensity $\A(\Ds\sqrt\dv)$ has the rank $k=4$,
 because the semidensity $\A(\sqrt\dv)$
 which is equal to (3.1) has the rank $k=2$.
 This is hidden in representation (3.8), where the function
  $\rho(z)$ corresponding to the volume
 form in adjusted coordinates depends
 non-explicitly on 
 derivatives of surface parametrization $z(\zeta)$.
 For example  consider the surface
  $x^0=\Psi\xi^1\eta_1,\theta_0=0,x^1=\xi^1,\theta_1=\eta_1$
  in the normal special odd symplectic superspace $E^{2.2}$
  with volume form $\dv=|dx d\theta|$.
   ($\Psi$ is an odd constant.)
  In Darboux coordinates 
  $\t x^0=x^0-\Psi x^1\theta_1$,
  $\t\theta_0=\theta_0$,
  $\t x^1=x^1+\Psi\theta_0 x^1$,
  $\t\theta_1=\theta_1-\Psi\theta_0 \theta_1$
  adjusted to the surface $M$,
  the volume form is equal to
  $\dv=(1+2\Psi\t\theta_0)
  |d\t x d\t\theta|$,
  hence the semidensity (3.8) is equal to
  $\Psi|d\t x^1 d\t\theta_1|^{1/2}$.

One can construct densities using 
 $\Ds$-operator
on embedded surfaces provided with induced symplectic structure.
For example one can consider on every non-degenerated surface 
$M^{n-1.n-1}$, the semidensity 
$\widetilde{\Delta_M^{\#}}\A(M,\,\sqrt\dv)$,
where $\widetilde{\Delta^{\#}_M}$ is the $\Ds$ operator on the
surface $M$ w.r.t. to the induced symplectic structure.
But in this case
                         $$ 
\A(M,\,\Ds\sqrt\dv)=
-\widetilde{\Delta_M^{\#}}\A(M,\,\sqrt\dv)\,.
	                       \eqno (3.14)
	             $$
This can be immediately checked
  in Darboux coordinates (3.5) adjusted to the surface $M$.

Finally we consider a simple example of these constructions.

 {\bf Example 2}. Let $E^{3.3}$ be a superspace associated to 
 the $3$-dimensional space $E^3$, $E^{3.3}=ST^*E^3$.
 We consider on $E^3$ the differential form
    $$
    w=-dx^0\wedge dx^1\wedge dx^2+b_0 dx^0+b_1 dx^1+b_2 dx^2\,.
    $$	     
    To this differential form there corresponds
    the semidensity 
          $$
  {\bf s}=\tau^{\#}(w)=
  (1+b_0\theta_1\theta_2+b_1\theta_2\theta_0+
  b_2\theta_0\theta_1)
  |dx^0dx^1dx^2d\theta_0d\theta_1d\theta_2|^{1/2}\,.
            $$
Let $M$ be a surface in $E^{3.3}$ which is defined by equations
 $x^0=\theta_0=0$. $M$ is associated to the space $E^2$ with coordinates 
 $(x^1,x^2)$.  
 The value of the odd semidensity $\A({\bf s})$ on $M$ is equal
 to $(b_2\theta_1-b_1\theta_2)|dx^1dx^2d\theta_1d\theta_2|^{1/2}$.
 This corresponds to differential form $b_1dx^1+b_2dx^2$--the pull-back
 of $w$ on $E^2$.
 The value of even semidensity $\A(\Ds{\bf s})$ on $M$ is equal
 to $(\p_2b_1-\p_1 b_2)|dx^1dx^2d\theta_1d\theta_2|^{1/2}$.
 This corresponds to differential form $d(b_1dx^1+b_2dx^2)$
 $=(\p_2b_1-\p_1 b_2)dx^1\wedge dx^2$--- the pull-back
 of $dw$ on $E^2$.	    
 
 The even density (volume form) on $M$ is equal to
  $P_0=(\p_2b_1-\p_1 b_2)^2|dx^1dx^2d\theta_1d\theta_2|$
 and odd density $P_1$ 
 $=(\p_2b_1-\p_1 b_2)(b_1\theta_2-b_2\theta_1)
 |dx^1dx^2d\theta_1d\theta_2|$.
\bigskip
\centerline {\bf Discussions}
\medskip

 We hope that considerations presented in this paper
 can be generalized for constructing densities
 in an odd symplectic superspace in higher order derivatives
 on surfaces of arbitrary  dimension 
 and for finding the complete set of local invariants of this
 geometry. In particular from considerations of Lemma follows
 that one can try to find non-trivial invariant densities
  on non-degenerated surfaces of codimension $(p.p)$ only if their
  rank is greater than $p$.
 It is interesting to analyze from this point 
 of view relations between
 geometrical interpretations presented in this paper and
 in [1] where some relations 
 of semidensity (3.1)
 with mean curvature in Riemanian geometry
  were indicated.
 
 The densities presented in formulae (3.13)
 are needed to be investigated more in details. Particularly
  one have to present explicit formulae for them and
  consider  
  the corresponding functionals over surfaces.
  These functionals are equal to zero in the case
   if the volume form in the ambient special odd
   symplectic superspace obeys to BV-master equation.
 Do Euler-Lagrange motion equations
 for these functionals equal identically to zero,
 as for usual Poincare-Cartan integral invariants
(1.3)?

  It can be interesting to note also
  that symmetry transformations of these functionals
  are not exhausted only by transformations induced by
  diffeomorphisms (2.15) of underlying space. General canonical
  transformations of superspace induce
  mixing of corresponding differential
  forms. (See considerations after formula (2.19).)

                        \bigskip
		      \centerline{\bf 4. Acknowledgment}
		      \medskip
 This work is highly stimulated by
very illuminating discussions with S.P. Novikov during
my talk on his seminar in Moscow. I am deeply grateful to him.

I want to express also my indebtness 
to I.A.Batalin, V.M.Buchstaber, A.S.Schwarz and I.V.Tyutin
for encouraging me to do this work.

I am deeply grateful for hospitality
and support
of the Abdus Salam International Centre for Theoretical
Physics where the essential part of this work was performed.

           \bigskip

           \centerline{\bf References}

                 \medskip
 [\CMP] O.M.Khudaverdian,  Odd Invariant Semidensity and
 Divergence-like Operators in an Odd Symplectic Superspace,
 {\it Comm. Math.Phys}., 1998, {\bf 198}, pp.591--606

[\BVa] I.A.Batalin, G.A.Vilkovisky:---
         Gauge algebra and Quantization.
      {\it Phys.Lett.}, 1981, {\bf 102B} pp.27--31.

[\BVb] I.A.Batalin, G.A.Vilkovisky:---
     Closure of the Gauge Algebra,
 Generalized Lie Equations and Feynman Rules.
  {\it Nucl.Phys.} {\bf B234}, 106--124, (1984).

[\JMPa] O.M. Khudaverdian---
 Geometry of Superspace with Even and Odd Brackets--
    {\it J.Math.Phys}., 1991, {\bf 32}, pp. 1934--1937.
  (Preprint of Geneve Universite. UGVA--DPT 1989/05--613).

 [\MPL] O.M.Khudaverdian, A.P.Nersessian---
   On Ge\-omet\-ry of Ba\-ta\-lin-
   Vil\-ko\-vis\-ky 
   For\-ma\-lism {\it Mod.Phys.Lett}., 1993, {\bf A8},
     No.25 pp.2377--2385.

  [\JMPd] O.M.Khudaverdian, A.P.Nersessian---
  Batalin--Vilkovisky Formalism and
  Integration Theory on Manifolds.
  {\it J. Math. Phys}., 1996,  {\bf37},pp.3713-3724.

[\SchCMPa] A.S.Schwarz---Geometry of Batalin--Vilkovisky 
Formalism.
 {\it Commun. Math. Phys.} ,1993,  {\bf 155},  pp.249--260.

[\Gayduk] A.V. Gayduk, O.M.Khudaverdian, A.S Schwarz:
  Integration on
  Sur\-faces in  Superspace.
  {\it Teor. Mat. Fiz.} {\bf 52}, 375--383, (1982).

 [\But] O.M. Khudaverdian, R.L.Mkrtchian---
   Integral Invariants of Buttin
  Bracket.
  {\it Lett. Math. Phys.}, 1989, {\bf 18}, pp. 229--234.

   [\Poin] O.M. Khudaverdian, A.S. Schwarz, Yu. S. Tyupkin---
Integral Invariants for Supercanonical Transformations.---
 {\it Lett. Math. Phys.}, 1981, {\bf 5}, pp. 517--522.

 [\BarScha] M.A.Baranov, A.S.Schwarz--- 
  Characteristic Classes of
 Supergauge Fields. 
 
 \noindent{\it Funkts. Analiz i ego pril.},
 {\bf 18}, No.2,  53--54, (1984).

 [\BarSchb] M.A. Baranov, A.S.Schwarz---
  Cohomologies of Supermanifolds.
 {\it Funkts. analiz i ego pril.}.
 {\bf 18}, No.3, 69--70, (1984).

 [\leites]. D.A. Leites--- {\it The Theory of Supermanifolds.}
   Karelskij Filial AN SSSR (1983).

 [\SchCMPb] A.S.Schwarz---  Supergravity, Complex Geometry and
 G-structures.
   {\it Commun. Math. Phys.}, {\bf 87}, 37--63, (1982).

[\Berezin] F.A.Berezin
{\it Introduction to Algebra and Analysis with Anticommuting
 Va\-ri\-ables.}  Moscow, MGU (1983).
 ({\it in English}--- Introduction to Superanalysis.
  Dordrecht--Boston: D.Reidel Pub. Co., (1987)).

[\leitb]. D.A.Leites: The new Lie Superalgebras and Mechanics.
   {\it Docl. Acad. Nauk SSSR} {\bf 236}, 804--807, (1977).

[\Shander] V.N.Shander---
  Analogues of the Frobenius and Darboux Theorems for Supermanifolds.
  {\it Comptes rendus de l'
 Academie bulgare des Sciences}, {\bf 36}, n.3,
 309--311, (1983).

 [\buttin] C.Buttin---C.R.Acad. Sci. Paris, Ser.A--B,
  1969, {\bf 269} A--87.

[\Gui]  V.Guillemin, S.Sternberg---{\it Geometric Asymptotics.}
AMS, Providence, Rhode Island, 1977. 
 (Mathematical Surveys. Number 14)

  \bye